\documentclass[11pt,dvipdfmx]{article}

\usepackage{amssymb}
\usepackage{amsthm}
\usepackage{amsmath}
\usepackage{ascmac}
\usepackage{color}
\usepackage{graphicx}
\usepackage{enumerate}
\usepackage{url}
\usepackage{caption} 

\usepackage[
bookmarks=true,
bookmarksnumbered=true, 
bookmarkstype=toc,
setpagesize=false,
pdftitle={A Successive LP Approach with C-VaR Type Constraints for IMRT Optimization},
pdfauthor={Makoto Yamashita},
pdfkeywords={Intensity-modulated radiotherapy treatment, Fluence map optimization, Linear Programming, C-VaR
}
]{hyperref}

\def\@themcountersep{}
\newtheorem{THEO}{Theorem}[section]
\newtheorem{ALGO}[THEO]{Algorithm}

\newtheorem{LEMM}[THEO]{Lemma}

\def\0{\mbox{\bf 0}}
\def\1{\mbox{\bf 1}}
\def\2{\mbox{\bf 2}}
\def\3{\mbox{\bf 3}}
\def\4{\mbox{\bf 4}}
\def\5{\mbox{\bf 5}}
\def\6{\mbox{\bf 6}}
\def\7{\mbox{\bf 7}}
\def\8{\mbox{\bf 8}}
\def\9{\mbox{\bf 9}}

\newdimen\zhige \zhige=0pt
\def\chige#1{{\setbox\zhige\hbox{#1}\ifdim\ht\zhige=1ex\accent24 #1%
  \else\ooalign{\unhbox\zhige\crcr\hidewidth\char24\hidewidth}\fi}}

\def\x{\mbox{\boldmath $x$}}

\def\z{\mbox{\boldmath $z$}}

\def\D{\mbox{\boldmath $D$}}

\def\FC{\mbox{$\cal F$}}

\def\SC{\mbox{$\cal S$}}

\def\Real{\mbox{$\mathbb{R}$}}

\ifx11 

\fi

\begin{document}

\noindent \textcolor{blue}{\leaders\hrule width 0pt height 0.08cm \hfill}
\vspace{0.5cm} 

{\Large \bf \noindent A Successive LP Approach with C-VaR Type Constraints for IMRT Optimization}

\vspace{0.2cm}
 \noindent
Shogo Kishimoto
\footnote{
Department of Mathematical and Computing Sciences,
 Tokyo Institute of Technology, 2-12-1-W8-29 Ookayama, Meguro-ku, Tokyo
 152-8552, Japan (kishimoto.s.ac@m.titech.ac.jp).
 } and
Makoto Yamashita
\footnote{
Department of Mathematical and Computing Science,
 Tokyo Institute of Technology, 2-12-1-W8-29 Ookayama, Meguro-ku, Tokyo
 152-8552, Japan (Makoto.Yamashita@is.titech.ac.jp).
His research was partially supported by JSPS KAKENHI (Grant Number: 15K00032).
 }
\\
Submitted: December 5, 2016.

\vspace{0.3cm}

\noindent {\bf Abstract:} 
Radiation therapy is considered to be one of important 
treatment protocols for cancers. 
Radiation therapy employs several beams of ionizing radiation to kill cancer tumors, but
such irradiation also causes damage to normal tissues.
Therefore, a treatment plan should satisfy dose-volume constraints (DVCs). 
Intensity-modulated radiotherapy treatment (IMRT) enables to control the beam intensities
and gives more flexibility for the treatment plan to satisfy the DVCs.
Romeijn et al. [Physics in Medicine and Biology, 48(21):3521, 2003] replaced the DVCs in an IMRT optimization with C-VaR (Conditional Value-at-Risk) type constraints, and proposed a numerical method based
on linear programming (LP). 
Their approach reduced the computation cost of the original DVCs, but
the feasible region of their LP problems was much narrow compared to the DVCs,
therefore, their approach 
often failed to find a feasible plan even when the DVCs were not so tight.

In this paper,  we propose a successive LP approach with the C-VaR type constraints.
We detect outliers form the solution of LP problems, and remove them from the domain of the C-VaR type constraints.
This eases the sensitivity of C-VaR type constraints to outliers and we can search feasible plans from wider regions.
Furthermore, we can give a mathematical proof that if the optimal value of the LP problem in the proposed approach is non-positive,
the corresponding  optimal solution satisfies all the DVCs.
From numerical experiments on test data sets, we observed that our approach found feasible solutions
more appropriately than existing LP approaches. In addition, our approach required fewer LP problems,
and this led to a short computation time.

\vspace{0.2cm}

\noindent \textcolor{blue}{\leaders\hrule width 0pt height 0.08cm \hfill}
\vspace{0.5cm} 
%
%
%
%
%

\section{Introduction}

In many countries, cancer is considered to be one of the principal causes of death. 
In Japan, it was reported in~\cite{MATSUDA13} that the fatalities number rose to 350 thousand people
and  800 thousand people were
newly diagnosed as cancer in the year 2010.
Prevalent types of cancer treatment include
chemotherapy, surgery, and radiation therapy.
An investigation conducted by Ministry of Health, Labor and Welfare of Japan
 \cite{mhlw2016} reported that
their percentages are 81\%, 72\%, and 32\%, respectively
(the numbers include combinations of treatment types).
National Cancer Institute also reported that half of the cancer patients 
receive radiation therapy during their treatment~\cite{nci2010}.
Radiation therapy is 
a treatment that uses several beams of ionizing radiation 
against cancer tumors, using a property that 
the beam irradiation has effect of reducing or killing cancer tumors. 
One of the merits of radiation therapy is that patients receive weak damage to the body
compared to surgery, and this can bring
high possibilities for continuing a normal daily life.

\begin{figure}[tbp]
\begin{center}
\includegraphics[scale=0.6]{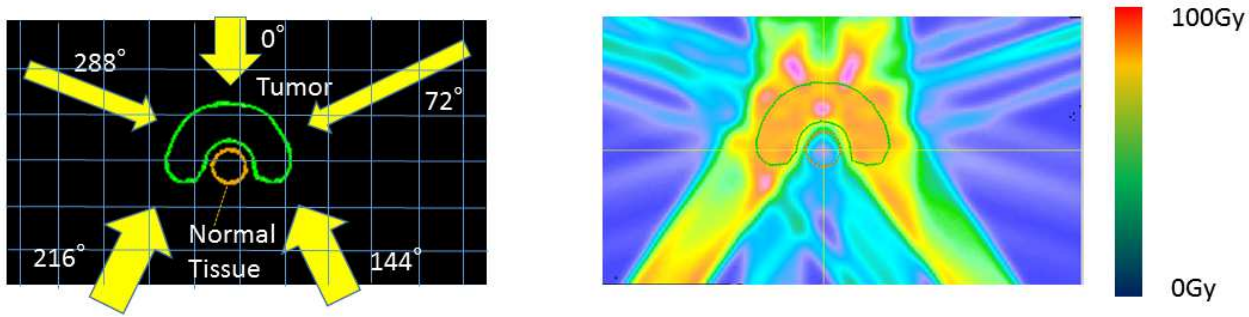}
\caption{An image of beam intensities for a tumor and normal tissues.} \label{fig:beam}
\end{center}
\end{figure}

Intensity-modulated radiotherapy treatment  (IMRT) has brought 
a remarkable flexibility in dose irradiated from the beams.
With the aid of computers 
and the usage of devices like multi-collimator, 
IMRT can control  the beam irradiation with higher accuracy than before,
and this leads to the high accuracy of the radiation delivery.
The irradiation device can continuously rotate around a patient who is on a couch 
and it irradiates the ionizing radiation 
at certain angles. In the left figure of Figure~\ref{fig:beam}, the beam is irradiated at the five angles,
$0^\circ, 72^\circ, 144^\circ, 216^\circ$ and $288^\circ$.
The right figure of Figure~\ref{fig:beam} shows that, due to an appropriate adjustment of the beam intensities,
the tumor receives a high dose and at the same time the normal tissues are kept away from the high dose.
The patient, however,  can not move during the treatment to increase the irradiation accuracy.
In practice, the number of angles during one treatment is limited to four to nine
in order to lessen a burden on the patient~\cite{dias2014genetic}.

The computation of IMRT planning involves several optimization aspects,
for example, 
beam angle optimization (BAO)~\cite{bertsimas2013hybrid, craft2007local, djajaputra2003algorithm, wang2004development, gao2016machinelearning},
fluence map optimization (FMO)~\cite{bednarz2002use, romeijn2003novel, tuncel2012strong},
and direct aperture optimization (DAO)~\cite{shepard2002direct, van2006intensity}.
BAO chooses the best angles from candidate angles on constraints like 
the number of angles available for one treatment, and
FMO is an optimization problem to determine the irradiation intensity for given beam angles. 
DAO considers the locations of multi-collimator and beam intensity simultaneously.
Chapters 4 and 5 of the handbook~\cite{pardalos2009handbook} discussed
many aspects on IMRT from the viewpoints of optimization.

A difficulty arising in FMO problems, however,  is that 
not only malignant tumors but also normal tissues near the tumors
receive negative effects from the beam irradiation.
The oncologists develop treatment plans for the irradiation areas and the beam intensity
to reduce the damage onto the normal tissues.

\begin{figure}[htbp]
\begin{center}
    \includegraphics[width=10.0cm]{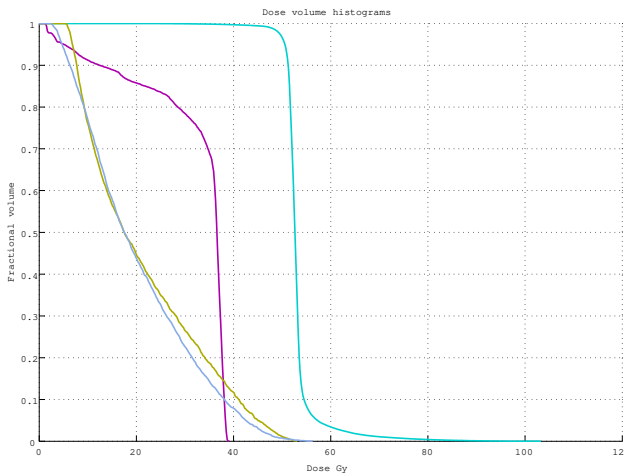}
\newline
\caption{An example of dose-volume histograms for a tumor (PTV) and normal tissues
(Cord, Lt Parotid and Rt Parotid).
} \label{fig:DVH}　
\end{center}
\end{figure}

A key criterion of the treatment plans is to satisfy dose-volume constraints (DVCs). 
A DVC is a radiation-dose constraint on a partial volume of an organ.
For example, 
more than 90\% area of the tumor should receives at least 
50 Gy (Gy is the international unit of ionizing radiation dose per unit mass; 1 Gy = 1 Joule/kg),
while 
the fraction of the normal tissues that receive 25 Gy or higher should be under 10\%.
Figure~\ref{fig:DVH} shows an example of dose-volume histograms (DVHs) for a cancer tumor 
and normal tissues near the tumor.
The blue line is a histograms for a tumor PTV and 
the other lines are for Cord, Lt Parotid and Rt Parotid.
In Figure~\ref{fig:DVH}, the horizontal axis is a dose volume and the vertical axis is 
the fraction of the structure.
The PTV histogram
passes the point (50 Gy, 90\%), and this indicates that the tumor area that receives 50 Gy or more is 90\%.
This histogram satisfies a DVC on the tumor which requires that the area
that incurs at least 50 Gy takes at least 90\% of the tumor (This DVC will be expressed 
as $L_{\textrm{PTV}}^{0.90} = 50$ in the notation that will be introduced in Section~2.1).

In general, multiple DVCs can be imposed for one treatment.
To find a good treatment plan that satisfies all the DVCs on tumors and normal tissues,
a number of approaches based on mathematical optimization methods have been proposed.
Morrill et al.~\cite{morrill1991dose} employed Linear programming (LP) problems for FMO problems with DVCs. 
Merrit et al.~\cite{merritt2002successive} proposed a successive LP approach.
In determining the beam intensity, a few parts of body receive much large dose, and 
such outliers make it very hard to find a feasible solution.
Merrit et al. detected these outliers based on the information of a dual LP problem,
and relaxed the dose threshold for the outliers to higher values.
Zhang et al.~\cite{zhang2010two} proposed a two-stage sequential LP approach.
Aleman et al.~\cite{aleman2010interior} solved an optimization problem that minimizes  a quadratic  objective function 
which evaluated the deviations from DVCs.
They proposed a specific interior-point point method to solve this quadratic optimization problem.
Hamacher and K{\"u}der~\cite{hamacher2002inverse} examined a multiple objective optimization approach .
Chu et al. \cite{chu2005robust} and Olafsson and Wright \cite{olafsson2006efficient} discussed approaches based on robust optimization framework with second-order cone programming problems.
Romeijn et al.~\cite{romeijn2003novel} 
introduced a concept of C-VaR (conditional value-at-risk) that had been originally developed in financial engineering.
Instead of DVCs, they used constraints that 
the \textit{average} dose for a given fraction part should satisfy a given threshold. 
As pointed in \cite{wu2002multiple},  FMO problems with DVCs formulated in mathematical optimization problems
often have multiple local minimum solutions
and such problems are essentially NP-hard~\cite{tuncel2012strong}.

Among the above methods, an advantage of the C-VaR type constraints of Romeijn et al.~\cite{romeijn2003novel} 
is that any irradiation intensity obtained from their LP problems satisfies all the DVCs. This property is not clearly mentioned 
in their paper, and we will verify it later in Lemma~\ref{le:C-VaR}. 
They solved only one LP problem, therefore, its computation cost was not expensive.
However, the feasible region of their LP problem was much 
narrow compared to the region intended by the DVCs. 
Since LP problems can be solved by a polynomial-time algorithms while the FMO problems themselves are NP-hard,
we can not completely remove the gap between the C-VaR type constraints and the DVCs.
In particular, the outliers seriously affected the average dose. 
This approach failed to find the beam intensities 
for some test instances of 
 Task Group (TG) 119 report by the American Association of Physicists in 
Medicine (AAPM)~\cite{ezzell2009imrt}.

In this paper, we propose a successive LP approach that employs the C-VaR type constraints.
We detect the outliers using LP problems which evaluate the deviation from the DVCs in the constraints,
and remove the outliers successively from the domain of the C-VaR type constraints.
We will show that if the objective value of the successive LP problems becomes non-positive,
the proposed method outputs beam intensities that satisfy all the DVCs.
The adjustment of the outliers by the successive LP problems enables  the proposed method to 
search a wider region than the approach of Romeijn et al.
In addition, the sequence of the objective values of the successive LP problems 
are non-increasing. This property implies that we can generate a sequence of the solutions that approaches to
the DVCs.
Since our optimization problems are still LP problems, the computation cost 
does not increase so much compared to the successive LP approach of Merrit et al.~\cite{merritt2002successive}.

We conducted numerical tests to verify the performance of our approach.
For test instances included in TG 119.  
our approach successfully found beam intensities that satisfied the DVCs
within a few LP iterations.
In addition, the solution of our approach  satisfied the DVCs more appropriately than Merrit et al.

The rest of this paper is organized as follow.
In Section~2, we first give a more precise definition and notation on DVCs
and briefly discuss the formulations of Merrit at al. and Romeijn et al.
In Section~3, we describe the details on our approach 
and, in Theorem~\ref{th:Update}, we give a proof of mathematical properties of our approach 
that are favorable for the FMO computation.
Section~4 reports the numerical results on the TG 119 test instances.
In Section~5, we will discuss several aspects of our approaches and extensions. 
Finally, we will give a conclusion in Section~6.

\section{Preliminaries and Existing Formulation}

As notation, we use $|S|$ to denote the cardinality of a set $S$.
We take a nonnegative part of a number $x$ by denoting $(x)^+ := \max\{x,0\}$.

\subsection{Preliminaries}

 To apply numerical computation to the IMRT problem in a practical way, 
the intended organs (or structures) and the radiation beams are often discretized into \textit{voxels} (small cuboid units)
 and \textit{beamlets}, respectively.
Let $S$ be the set of the structures.
For each structure $s \in S$, we use $V_s$ to denote the voxel set of $s$.
Without loss of generality, we assume $|V_s| > 0$ throughout of this paper.  
In IMRT, the multi-leaf collimators make it possible to treat the radiation beam as a set of beamlet
and control the intensity of each beamlet independently. The set of beamlet is denoted by $B$.
The area radiated from each beamlet is usually 10 mm $\times$ 10 mm,
5 mm $\times$ 5 mm, 3 mm $\times$ 3 mm
and the number of beamlets $|B|$ is from hundreds to thousands.
The order of the number of voxels is usually $10^4$.

It is often assumed that a dose of each voxel received from the beamlets is  a linear function,
therefore,  the dose that the $i$th voxel in the $s$th structure receives can be expressed in
$z_{si} = \sum_{j \in B} D_{sij} x_j$. Here, $x_j$ is the intensity of the $j$th beamlet.
The element $D_{sij}$ is the $(i,j)$th element of a matrix 
$\D_s \in \Real^{|V_s|\times |B|}$, and the matrix $\D_s$ is called a fluence matrix and its element $D_{sij}$ expresses
the dose that the $i$th voxel in the $s$th structure
receives from the unit intensity of the $j$th beamlet.
To compute the fluence matrix, Naqvi et al.~\cite{naqvi2003convolution}
utilized the Monte Carlo superposition.
In the material below, we assume that the fluence matrix is given.
 
In a suitable treatment plan, its corresponding histogram should satisfy the DVCs.
Precisely speaking, a DVC is identified by a structure $s$ and
a fractional parameter $\alpha \in (0,1)$. 
We can classify the DVCs into the two types, the upper and lower DVCs;\newline
\begin{tabular}{lp{14cm}}
 (upper) & The fraction of the voxels in the structure $s$ that receive at least $U^{\alpha}_s$ Gy is at most $\alpha$. \\
 (lower) & The fraction of the voxels in the structure $s$ that receive at least $L^{\alpha}_s $ Gy is at least $\alpha$.
\end{tabular}
\newline
As an example, let us impose three DVCs to a tumor named PTV;
$L_{\text{PTV}}^{0.9}  =  50.0$, $L_{\text{PTV}}^{0.99} =  46.5$
and $U_{\text{PTV}}^{0.2} = 55.0$.  
In this case, at least $99\%$ of the tumor PTV must receive at least 
$46.5$ Gy. Furthermore,  $90\%$ of  PTV should receive a higher dose than 
 $50.0$ Gy. At the same time, we should also avoid extremely strong intensity 
and this is expressed by the upper
constraint of $U_{\text{PTV}}^{0.2}$, that is, 
 $20\%$ voxels or less of PTV can exceed $55.0$ Gy.
As shown in this example, the number of DVCs
imposed for one tumor or organ can be larger than one.

In the following discussion, we will use $\underline{A}_s$ and $\overline{A}_s$ 
to denote the set of the fractions that are involved in the lower and upper DVCs of the structure $s$, respectively.
For each $\alpha \in \overline{A}_s$, we associate the upper DVC whose threshold is $U_s^{\alpha}$ Gy.
Similar notation is applied to $ \alpha \in \underline{A}_s$ for the lower DVC with $L_s^{\alpha}$ Gy.
In the DVC example above, we have $\underline{A}_{PTV} = \{0.9, 0.99\}$ and $\overline{A}_{PTV} = \{0.2\}$.
In a mathematical form, a single upper DVC can be described as 
$\frac{\left| \{ i \in V_s | z_{si} > U_s^{\alpha} \} \right|}{|V_s|} \le \alpha$,
and a single lower DVC as 
$\frac{\left| \{ i \in V_s | z_{si} > L_s^{\alpha} \} \right|}{|V_s|} \ge \alpha$.

We assume that $\underline{A}_s \subset (0,1) $ and $\overline{A}_s \subset (0,1)$.
For the specific fractional case corresponding to $\alpha = 0$ or $\alpha = 1$, 
we also denote the upper or lower bounds by
$U_s$ and $L_s$, respectively. 
When these threshold are used, 
\textit{each} voxel in the structure $s$ is required to receive the dose between $L_s$ and $U_s$.
 
 
%
%
%
%
　
An FMO problem can now be casted as a mathematical problem to find
the beamlet intensities that satisfy all the DVCs.
If we are allowed to use mixed-integer programming problems, 
one goal in FMO is to find a solution 
of the feasible set $\FC$ defined by
\begin{eqnarray}\label{eq:FMO}
\begin{array}{llcll}
\FC := & \{\x \in \Real^{|B|} &:& 
   \sum_{j = 1}^{|B|} D_{sij} x_j = z_{si} & \mbox{for }  i = 1,\ldots, |V_s|; s = 1,\dots, |S|\\
   & & & L_s \leq z_{si} \leq U_s &  \mbox{for } i = 1,\ldots, |V_s|; s = 1, \ldots, |S| \\
   & & & z_{si} \geq 0 &  \mbox{for } i = 1,\ldots, |V_s|; s = 1, \ldots, |S| \\
   & & & z_{si} \ge L_s^{\alpha} \underline{b}_{si}^{\alpha} 
 &  \mbox{for } i = 1,\ldots, |V_s|; \alpha \in \underline{A}_s; s = 1, \ldots, |S| \\
   & & & \underline{b}_{si}^{\alpha} \in \{0,1\} &  \mbox{for } i = 1,\ldots, |V_s|; \alpha \in \underline{A}_s; s = 1, \ldots, |S|\\
   & & & \sum_{i=1}^{|V_s|} \underline{b}_{si}^{\alpha} \ge \alpha |V_s| 
 & \mbox{for }  \alpha \in \underline{A}_s; s = 1, \ldots, |S| \\
   & & & z_{si} \le U_s^{\alpha}  + M \overline{b}_{si}^{\alpha} 
 &  \mbox{for } i = 1,\ldots, |V_s|; \alpha \in \overline{A}_s; s = 1, \ldots, |S| \\
   & & & \overline{b}_{si}^{\alpha} \in \{0,1\} &  \mbox{for } i = 1,\ldots, |V_s|; \alpha \in \overline{A}_s; s = 1, \ldots, |S| \\
&    & & \sum_{i=1}^{|V_s|} \overline{b}_{si}^{\alpha} \le \alpha |V_s| 
 &  \mbox{for } \alpha \in \overline{A}_s; s = 1, \ldots, |S| \\
 & & & x_j \ge 0 &  \mbox{for } j = 1, \ldots, |B| \quad \}.
\end{array}
\end{eqnarray}

In this definition, 
a scalar $M$ is so-called big-M, a constant number large enough.
To express the fraction of the partial volume,
the binary variables 
$\underline{b}_{si}^{\alpha}$ and $\overline{b}_{si}^{\alpha}$ are introduced.
We should remark that a single upper DVC 
$|\left\{ i \in V_s | z_{si} > U_s^{\alpha} \right\}| \le \alpha |V_s|$
is imposed by a combination of 
$z_{si} \le U_s^{\alpha}  + M \overline{b}_{si}^{\alpha}$, 
$\overline{b}_{si}^{\alpha} \in \{0,1\}$ and 
$\sum_{i=1}^{|V_s|} \overline{b}_{si}^{\alpha} \le \alpha |V_s|$.
The number of voxels exceeds thousands in practical situations, so 
the number of these binary variables are also considerably large.
The set $\FC$ embraces properties of combinatorial sets 
and it is not an easy task to find a feasible point of $\FC$ exactly.
As pointed in~\cite{tuncel2012strong}, such a task is NP-hard.

\subsection{A successive linear programming method}
In 2002, Merritt et al.~\cite{merritt2002successive}
employed LP problems to formulate the FMO problem in a mathematical way
and exploited the information obtained from the dual LP problems.
Their method is refereed as Method-M in the material below.
A main idea of Method-M is to solve LPs successively searching better beamlet intensities.
 We now briefly introduce Method-M by a simple situation which involves
one tumor structure ($s = 1$) and one healthy structure ($s = 2$).
We consider hard DVC thresholds  $U_1$ and $U_2$ on the tumor and the healthy structures, respectively,
and use a soft DVC threshold $U_2^{\alpha}$ such that $U_2^{\alpha} \le U_2$.
A framework of Method-M for this situation is given as follow, and
this framework can be easily extended to general cases of more DVCs.

\begin{ALGO}\label{al:Merritt}~{\upshape\cite{merritt2002successive}}
(Method-M) 
A successive LP method for the FMO optimization
\end{ALGO}
\begin{enumerate}
\item Set the iteration number $k = 0$, and set an initial set $R_0 = \emptyset$. Set a parameter $\lambda > 0$
and a stopping threshold $\hat{\tau} > 0$.

\item Solve the following $k$th LP to determine the intensity of beamlets $\x \in \Real^{|B|}$
and let  $\tau_k$ be the optimal value of this LP.
\begin{align*}
  &&&&& \textrm{maximize}   && \tau && &&&& \\
  &&&&& \textrm{subject to} && \sum_{j=1}^{|B|} D_{sij} x_j = z_{si} && \mbox{ for } i=1,\dots, |V_s| ;  s = 1,2\\
  &&&&&                   && \tau \leq z_{1i} \leq U_1  && \mbox{ for } i=1,\dots, |V_1| &&&&& \\
  &&&&&                   && z_{2i} \leq U_2 && i \in R_k \\
  &&&&&                   && z_{2i} \leq U_2^{\alpha} && i \notin R_k \\
  &&&&&                   && z_{si}\geq 0 && \mbox{ for } i=1,\dots, |V_s|;  s = 1,2\\
  &&&&&                   && x_{j}\geq 0 && \mbox{ for } j = 1, \ldots, |B|.
\end{align*}

\item If $\tau_k > \hat{\tau}$, output the optimal solution of the $k$th LP and stop.
\item Update  $R_k$ with the rule
$R_{k+1} = R_k \cup \{i \in V_2 \  | \  y_i^* > \lambda \},$
where $y_i^*$ is the dual variable corresponding to the constraint $z_{2i} \leq U_2^{\alpha}$ in the $k$th LP.
\item Increment $k$ and return to Step~2.
\end{enumerate}

In the $k$th LP, each voxel in the tumor structure ($s = 1$) receives the dose 
at least $\tau$,
therefore, the aim in the $k$th LP is to maximize $\tau$ so that the voxels in the tumor receives as high dose as possible.

The usage of the hot-spot set $R_k$ characterizes
Method-M. Though it is preferable that each voxel in the healthy 
structure ($i \in V_2$) should be imposed by a low limit $z_{2i} \le U_2^\alpha$, 
such a constraint is too restrictive to satisfy.
Therefore, 
Merritt et al. relaxed this constraint so that a small set $R_k \subset V_2$ can be exposed to higher dose than $U_2$ 
(Note that $U_2 \ge U_2^{\alpha}$).
These voxels, regarded as outliers, are determined using the information from the dual LP problem in Step~4. 
The updated hot-spot set $R_{k+1}$ is composed of the voxels for which the lower limit $z_{2i} \le U_2^\alpha$
is too restrictive.

In contrast, a disadvantage of this method is that it does not take the fractional parameter $\alpha$ of tumor into consideration.
The constraints in the LP problem involve all the voxels in the structures, and they  
are much stronger than the DVCs. 
Since the feasible region has a tendency to become very narrow, this method may fail to find favorable beamlet intensities 
for DVCs.

\subsection{An Approach Based on C-VaR type Constraints}
Romeijn et al. \cite{romeijn2003novel, romeijn2006new} also utilized LP problems to determine the beamlet intensities,
but their approach brought a different perspective. 
Their method is refereed as Method-R in this paper.
The key step of Method-R is to replace the DVCs  with C-VaR type constraints in the LP problem,
and only one LP problem is solved.
C-VaR stands for \textit{conditional value at risk}, and it was originally introduced 
by Rockafellar \cite{rockafellar2000optimization} in a context of economics. 
In economics, there is a demand  for estimating the expected value
of the loss that exceeds a certain level called Value at Risk (VaR).
In particular, C-VaR has high affinity with a fraction.

For a random variable $X$ and level $\alpha$, 
the original definition of C-VaR computes an average of VaR using an integral as follow:
\[\text{CVaR}_{1-\alpha}(X) := \frac{1}{\alpha} \int_{0}^{\alpha} \text{VaR}_{1-\tau }(X) d\tau. \]

Here, we do not discuss a precise definition of VaR, since  an equivalent but  more convenient form of C-VaR is available:
\begin{eqnarray}
\text{CVaR}_{\alpha}(X) := \min_{C \in \mathbb{R}} \left\{C +\frac{1}{1-\alpha} E[(X-C)^+] \right\}. \label{eq:C-VaR}
\end{eqnarray}
A key step of Method-R 
is to replace a single upper DVC 
$\left| \left\{ i \in V_s | z_{si} > U_s^{\alpha} \right\} \right| \le \alpha |V_s|$
with a C-VaR type inequality
\begin{eqnarray*} 
\min_{C \in \Real}\left\{ C 
+ \frac{1}{\alpha |V_s|}\sum_{i=1}^{|V_s|} (z_{si} - C)^+\right\}
\le
U_s^{\alpha}.
\end{eqnarray*}
This inequality is equivalent to find $\overline{\zeta}_s^{\alpha} \in \Real$ which satisfies
\begin{eqnarray}
\overline{\zeta}_s^{\alpha}
+ \frac{1}{\alpha |V_s| }\sum_{i=1}^{|V_s|} (z_{si} - \overline{\zeta}_s^{\alpha})^+
\le
U_s^{\alpha}.\label{eq:upperDVC2}
\end{eqnarray}

The upper DVC imposes that the number of voxels that receive  $U_s^{\alpha}$ or more is bounded by $\alpha |V_s|$.
In contrast, the upper C-VaR type constraint~(\ref{eq:upperDVC2}) requires
that the mean dose of the $\alpha |V_s|$ highest voxels 
be under $U_s^{\alpha}$.
\begin{figure}
     \centering 
     \includegraphics[width=4.0cm]{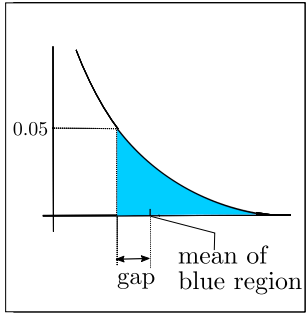}
 	   \caption{A comparison between DVC and C-VaR type constraint}
 	   \label{fig:CVaR}
\end{figure}
Here, we use Figure~\ref{fig:CVaR} to compare DVC and the C-VaR type constraint
in the DVH style. 
The upper DVC on $U_s^{\alpha} = 0.05$ requires the leftmost point of blue region be less than $U_s^{\alpha}$.
In contrast, the upper C-VaR corresponding to this DVC requires that the mean of blue region be less than $U_s^{\alpha}$.
Therefore, when $U_s^{\alpha}$ lies between the leftmost point and the mean,
there is a difference between the DVC and the C-VaR type constraint.
In Chapter~3, we will observe this in detail.

We can develop the inequality (\ref{eq:upperDVC2}) in another way. 
For simplicity, we assume that 
$\alpha |V_s|$ is an integer number.
Let us use $\SC_{k}(\z)$ for an integer $k$ and $\z \in \Real^n$
to denote the sum of the $k$ largest elements of $\z$.
In~\cite[Example 9.10]{calafiore2014optimization}, it is shown that $\SC_{k}(\z) \le t$ is equivalent to a condition that there exists
$\overline{\zeta}, t_1, t_2, \ldots, t_n \in \Real$  such that 
\begin{eqnarray}
\left\{\begin{array}{l}
t \ge  k \overline{\zeta} + (t_1 + t_2 + \cdots + t_n )  \\
t_1, t_2, \ldots, t_n \ge 0 \\
t_1 - z_1 + \overline{\zeta} \ge 0, t_2 - z_2 + \overline{\zeta} \ge 0, t_n - z_n + \overline{\zeta} \ge 0.
\end{array}\right. \label{eq:largests}
\end{eqnarray}
Since the mean dose of the $\alpha |V_s|$ highest voxels should be less than $U_s^{\alpha}$,
we have  $\frac{\SC_{\alpha |V_s|} (\z_s)}{\alpha |V_s|} \le U_s^{\alpha}$.
Here, $\z_s$ is the vector that collects the doses in the structure $s$,
$\z_s :=(z_{s1}, \ldots, z_{s |V_s|})^T \in \Real^{|V_s|}$.
Applying (\ref{eq:largests}) with $k = {\alpha |V_s|}$
and $t = {\alpha |V_s|} \cdot U_s^{\alpha}$, we can show that 
$\frac{\SC_{\alpha |V_s|} (\z_s)}
{\alpha |V_s|} \le U_s^{\alpha}$ if and only if 
there exists $\overline{\zeta}_s^{\alpha} \in \Real$ 
that satisfies (\ref{eq:upperDVC2}). 

By introducing  the concept of C-VaR to the IMRT optimization, 
Romeijn et al.~\cite{romeijn2003novel} proposed 
the following LP problem:
\begin{subequations}
\label{eq:RomeijnLP}
\begin{align}
  &&&&&\textrm{minimize}    && \sum_{s = 1}^{|S|} \sum_{i =1}^{|V_s|} F_s(z_{si}) & &&&&& \\
  &&&&& \textrm{subject to} && \sum_{j = 1}^{|B|} D_{sij} x_j = z_{si} & i = 1,\ldots, |V_s|; s = 1,\dots, |S|\\
  &&&&&                     && L_s \leq z_{si}\leq U_s &  i = 1,\ldots, |V_s|; s = 1, \ldots, |S| \label{hardDVC} \\
  &&&&&                     && \underline{\zeta}_{s}^{\alpha} - \frac{1}{(1-\alpha) |V_s|} \sum_{i = 1}^{|V_s|} ( \underline{\zeta}_{s}^{\alpha} - z_{si} ) ^{+}  \geq L_s^{\alpha}  & \alpha \in \underline{A}_s ; s = 1, \ldots,  |S| \label{lowerCVaR}\\
  &&&&&                     && \overline{\zeta}_{s}^{\alpha} + \frac{1}{\alpha |V_s|} \sum_{i=1}^{v_s} (z_{si} - \overline{\zeta}_{s}^{\alpha}  )^{+}  \leq U_s^{\alpha}  & \alpha \in \overline{A}_s ; s = 1, \ldots, |S| \label{upperCVaR}\\
  &&&&&                     && x_j  \geq 0   & j = 1, \ldots, |B| &&&&& \\
  &&&&&                     && z_{si}  \geq 0  & i = 1,\ldots, |V_s|; s=1, \ldots, |S|    &&&&& \\
  &&&&&                    &&  \underline{\zeta}_{s}^{\alpha},  \overline{ \zeta}_{s}^{\alpha} : \textrm{free variable} &&   
\end{align}
\end{subequations}
Here, the decision variables are  the beamlet intensities $x_1, \ldots, x_{|B|}$. 
To implement the C-VaR type constraints, intermediate variables
$\underline{\zeta}_{s}^{\alpha}$ and $\overline{ \zeta}_{s}^{\alpha}$ are employed. 

They used piecewise a linear function $F_s$ for the objective function in order
to express an deviation from their desired situation, and this objective function
remained (\ref{eq:RomeijnLP}) as an LP problem. 
Aleman et al.~\cite{aleman2010interior} examined quadratic penalty functions for the objective function
to incorporate the deviation,
and they applied interior-point methods to solve the resultant quadratic optimization problem.

The validity of C-VaR type constraints in (\ref{eq:RomeijnLP}) 
can be guaranteed by the following lemma.
Though this claim was partially implied in \cite{romeijn2006new}, we give it in an explicit way.

\begin{LEMM}\label{le:C-VaR}
Any feasible solution of (\ref{eq:RomeijnLP}) fulfills all the DVCs.
\end{LEMM}
\textbf{Proof:} From the constraint (\ref{hardDVC}), it is clear that 
the hard DVCs ($L_s \le z_{si} \le U_s$) are satisfied.

We now assume that a single upper DVC 
$|\left\{ i \in V_s | z_{si} > U_s^{\alpha} \right\}| \le \alpha |V_s|$
is violated, and we will derive a contradiction. 
From this assumption, the  number of voxels such that $ z_{si} > U_s^{\alpha}$ 
is greater than $\alpha |V_s|$.
Since $ U_{s}^{\alpha} \geq \overline{\zeta}_{s}^{\alpha}$ from 
(\ref{upperCVaR}), it holds that $z_{si} > \overline{\zeta}_{s}^{\alpha}$ when
$z_{si} > U_{s}^{\alpha}$.
There exists at least one $i$ such that $z_{si} > \overline{\zeta}_{s}^{\alpha}$, since $|\left\{ i \in V_s | z_{si} > U_s^{\alpha} \right\}| > \alpha |V_s| \ge 0$ 
from our assumption, hence,
 we have $(z_{si} - \overline{\zeta}_{s}^{\alpha})^+ 
= z_{si} - \overline{\zeta}_{s}^{\alpha} > 
 U_{s}^{\alpha} - \overline{\zeta}_{s}^{\alpha}$ for such $i$.
It leads to 
\begin{eqnarray*}
\overline{\zeta}_{s}^{\alpha} + \frac{1}{\alpha |V_s|} \sum_{i=1}^{|V_s|} (z_{si} - \overline{\zeta}_{s}^{\alpha}  )^{+}  
> 
 \overline{\zeta}_{s}^{\alpha} + \frac{1}{\alpha |V_s|} \left(\alpha |V_s|\right)  (U_s^{\alpha} - \overline{\zeta}_{s}^{\alpha}  )  
=  U_s^{\alpha} 
\end{eqnarray*}
and this contradicts to  (\ref{upperCVaR}). 
Hence, any feasible solution of (\ref{eq:RomeijnLP}) does not violate the upper DVC.
A similar discussion can be applied to the lower DVCs
 (\ref{lowerCVaR}).
\hfill \qed

This lemma indicates
an advantage of the LP model (\ref{eq:RomeijnLP}) that 
if a feasible solution is found in (\ref{eq:RomeijnLP}), it should fulfill all the DVCs.
A negative side is that this approach may fail to find solutions that exist in the gap between the DVCs and the C-VaR type constraints
(The converse of Lemma~\ref{le:C-VaR} does not hold in general).
This approach searches only narrower feasible region than the original DVCs,
and this aspect motivated us to develop a successive method that extends 
the feasible region from the C-VaR type constraints.

\section{Successive Linear Programming Approach with C-VaR type Constraints}

A main difficulty in handling the DVCs is that the DVCs involve the fraction $\alpha$, hence
finding a feasible solution is already a demanding task.
Method-M utilized the hot-spot to remove some outliers from the strict constraints,
and Method-R introduced the C-VaR type constraints to replace the original DVCs.
A common problem arising from the two methods 
was that their search regions were not wide enough to cover the region shaped by the DVCs.
Hence, the two methods sometimes fail to find a feasible solution,
even when the original feasible set 
has enough space.

In particular, we observed from preliminary numerical tests that the C-VaR type constraints are 
sensitive to outliers in the sense that 
the voxels that have extremely high or low doses affect the constraints of (\ref{eq:RomeijnLP}) seriously.

%

\begin{figure}[htbp]
     \centering 
     \centering
     \captionsetup{justification=centering}
 	   \includegraphics[width=4.0cm]{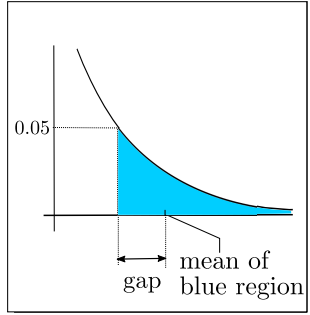}
 	   \caption{A DVH in which a few voxels receive extremely high doses }
 	   \label{fig:outliers_right}
\end{figure}


To illuminate such a phenomenon, we compare
Figures~\ref{fig:CVaR} and \ref{fig:outliers_right} that partially illustrate different solutions of an FMO problem in the style of DVH. 
The two DVHs are almost same, but Figure~\ref{fig:outliers_right} has a few voxels that receive extremely high doses.
%
Their means of the blue region are quite different, and it indicates that 
the satisfiability of (\ref{upperCVaR}) strongly depends on 
the voxels that have the highest doses.
Such voxels should be handled carefully as outliers.

We propose a method that combines the successive update of outliers 
and the concept of the C-VaR type constraints.
In the proposed method, 
we first solve an LP problem by relaxing the C-VaR constraints  so that this LP problem always has a feasible point.
From the optimal solution of this LP problem, we extract the outliers and transfer them from 
the domain of the C-VaR type constraints (\ref{upperCVaR})
to the sets of outliers.
Using the new sets, we build a next LP problem.
We repeat the updates of outliers until we obtain a feasible solution 
that satisfies the original DVCs. 
This new method has several favorable properties, and we will discuss them in Theorem~\ref{th:Update}.

\begin{figure}[htbp]
    \centering
    \includegraphics[width=8.0cm]{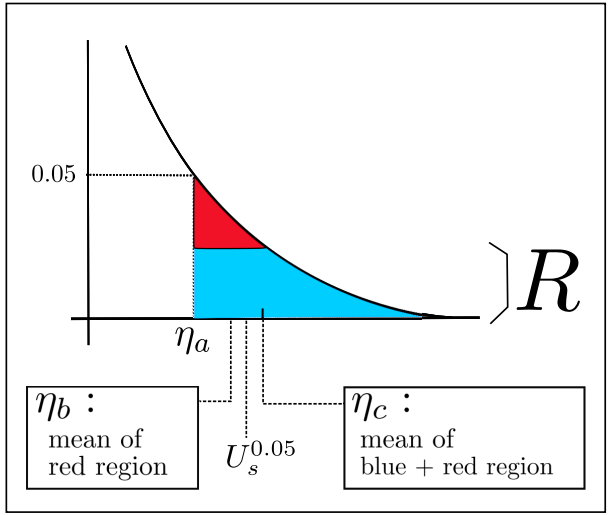}
    \caption{Effect of the exclusion of $R$}
    \label{fig:comparison_constraint}
\end{figure}

The comparison between an original DVC, its corresponding C-VaR type constraints and the constraint in the proposed
method is summarized as an illustrative example of~Figure~\ref{fig:comparison_constraint}.
We denote the left endpoint of the blue and red areas by $\eta_a$.
We also use $\eta_b$ and $\eta_c$ to denote the mean dose of the red area and that of blue and red areas, respectively.
The value on this DVH at $\eta_a$ is 0.05,
therefore this solution satisfies an DVC of $U_s^{0.05}$
if and only if $\eta_a \leq U_s^{0.05}$.
Similarly, $\eta_c \leq U_s^{0.05}$ if and only if
the solution satisfies the C-VaR type constraint for $U_s^{\alpha}$.
From the discussion in the previous paragraphs,
if $U_s^{0.05}$  in the interval $\eta_a < U_s^{0.05} < \eta_c$, 
this solution is a feasible solution of the FMO problem,
even though the approach based on the C-VaR type constraints fails to recognize 
this solution as a feasible solution.

We consider an effect of excluding voxels that have highest doses
by defining them as a set of outliers $R$.
We evaluate the left-hand side of the C-VaR type constraint using only the rest of the voxels.
In  Figure~\ref{fig:comparison_constraint}, this corresponds to the removal of the blue regions,
therefore, the left-hand side of the C-VaR type constraint is shifted from 
$\eta_c$ to $\eta_b$.
If $U_s^{0.05}$ is in the interval  $\eta_b < U_s^{0.05} < \eta_c$, 
this new approach can detect this solution is feasible.
This is the key idea of our approach.

The framework of the proposed method is summarized in Algorithm~\ref{al:proposed}.
In the $k$th LP problem of Step~2,
we use $\overline{R}_s^{k,\alpha}$  and $\underline{R}_s^{k,\alpha}$
to denote the sets of outliers
with respect to the thresholds $\alpha \in \overline{A}_s$ and 
$\alpha \in \underline{A}_s$, respectively.
In addition, the objective function $t$ is introduced to measure the deviation from the DVCs.
The positive constants $\overline{P}_s$, $\underline{P}_s$, $\overline{P}_s^{\alpha}$ 
 and $\underline{P}_s^{\alpha}$ are embedded to control the relaxation of  the DVCs.


\begin{ALGO}\label{al:proposed}
A successive updates of outliers with C-VaR type constraints for  FMO problems
\end{ALGO}
\begin{enumerate}
\item Set the iteration counter $k = 1$ and the initial sets of outliers
$\overline{R}_s^{1,\alpha} = \emptyset$ for $\alpha \in \overline{A}_s, s \in S$
and 
$\underline{R}_s^{1,\alpha} = \emptyset$ for $\alpha \in \underline{A}_s, s \in S$.
Choose positive constants $\overline{P}_s$ and $\underline{P}_s$ 
for $s \in S$, $\overline{P}_s^{\alpha}$  for $\alpha \in \overline{A}_s, s \in S$,
and $\underline{P}_s^{\alpha}$  for $\alpha \in \underline{A}_s, s \in S$.
\item Solve the following $k$th LP. Let $t^{(k)}$ be the optimal value of this LP problem, 
and $x_j^{(k)}$ and $z_{si}^{(k)}$ the obtained  solution.
\begin{subequations}\label{eq:NewProblem}
\begin{align} 
  &&&&&\textrm{min}    && t & &&&&& \\
  &&&&& \textrm{s.t.} && \sum_{j=1}^{|B|} D_{sij} x_j = z_{si} & i \in V_s;  s \in S \\
  &&&&&                     && z_{si}\leq U_s + \overline{P}_s t &  i \in V_s; s \in S \label{eq:6c} \\
  &&&&&                     && z_{si}\geq L_s - \underline{P}_s t & i \in V_s; s \in S \\
    &&&&&                     &&  \underline{\zeta}_{s}^{\alpha} -  \frac{1}{(1-\alpha) |V_s| 
 - |\underline{R}_s^{k,\alpha}|}  \sum_{\substack{i=1 \\ i \notin \underline{R}_s^{k,\alpha} }}^{|V_s| } (\underline{\zeta}_{s}^{\alpha} - z_{si} )^{+} 
   \geq L_s^{\alpha} - \underline{P}_s^{\alpha}t  &  \alpha \in \underline{A}_s; s \in S \label{New_lower_C} \\
  &&&&&                     &&  \overline{\zeta}_{s}^{\alpha} + \frac{1}{\alpha |V_s| 
 - |\overline{R}_s^{k,\alpha}|} \sum_{i=1,  i \notin \overline{R}_s^{k,\alpha}}^{|V_s|} (z_{si} - \overline{\zeta}_{s}^{\alpha}  )^{+} 
   \leq U_s^{\alpha} + \overline{P}_s^{\alpha}t  &  \alpha \in \overline{A}_s; s \in S  \label{New_upper_C} \\
  &&&&&                     && x_j  \geq 0   & j \in B &&&&& \\
  &&&&&                     && z_{si}  \geq 0  & i \in V_s; s \in S    &&&&& \\
  &&&&&                    &&  
 \overline{\zeta}_{s}^{\alpha}: \textrm{free variable} &
 \alpha \in \overline{A}_s; s \in S  \label{New_upper_zeta} &&&& \\
  &&&&&                    &&  
 \underline{\zeta}_{s}^{\alpha} : \textrm{free variable} &
 \alpha \in \underline{A}_s; s \in S  \label{New_lower_zeta} &&&&  \\
  &&&&&                    &&  
 t : \textrm{free variable} & &&&& 
\nonumber
\end{align}
\end{subequations}
\item If $t^{(k)} \le 0$, output $\x^{(k)}$ as the solution and stop.
\item 
Update the sets of outliers by the rules
\begin{eqnarray*}
\overline{R}_s^{k+1, \alpha} := \left\{ i \in V_s : z_{si}^{(k)} > U_s^{\alpha} + \overline{P}_s^{\alpha} t^{(k)} \right\},
\quad 
\underline{R}_s^{k+1, \alpha} := \left\{ i \in V_s : z_{si}^{(k)} < L_s^{\alpha} - \underline{P}_s^{\alpha} t^{(k)} \right\}.
\end{eqnarray*}
Increment $k$ and return to Step 2.
\end{enumerate}


The proposed method has the following suitable properties.
\begin{THEO}\label{th:Update}
We assume that the feasible region $\FC$ of the FMO problem (\ref{eq:FMO}) is not empty.
Then, for the LP problems solved in Algorithm~\ref{al:proposed}, it holds that 
\begin{enumerate}[(a)]
\item 
For any $k \ge 1$, the $k$th LP problem (\ref{eq:NewProblem}) has an optimal solution.
\item If $t^{(k)} \le 0$, the output solution $\x^{(k)}$ satisfies all the DVCs (that is, $\x^{(k)} \in {\FC}$).
\item The sequence $\{t^{(k)}\}$ is monotone non-increasing.
\end{enumerate}
\end{THEO}

Part (a) indicates that the solution $\x^{(k)}$ is well-defined through the execution of Algorithm~\ref{al:proposed}. 
Part (b) gives a validity for the stopping criterion $t^{(k)} \le 0$ in Step 3.
Finally, we can infer from Part (c)
that the sequence $\{\x^{(k)}\}$ has a tendency to approach to the set that satisfy all the DVCs.

We remark that 
the solution obtained from Method-R corresponds to $x^{(1)}$ of Algorithm~\ref{al:proposed}
with the parameters $\overline{P}_s = \underline{P}_s = \overline{P}_s^{\alpha} = \underline{P}_s^{\alpha} = 0$ 
for all $s$ and $\alpha$.
From Part(c), therefore, the proposed method is more flexible than Method-R.

\noindent \textbf{Proof:} 
For Part(a), we first consider the case $k = 1$ and discuss $k \ge 2$ by induction.
At the beginning of $k=1$, $\underline{R}_s^{1, \alpha}$ and $\overline{R}_s^{1, \alpha}$ are empty sets. 
The denominators, therefore,  in (\ref{New_lower_C}) and (\ref{New_upper_C}) 
are not zero,
as $|V_s| > 0$ without loss of generality and $0 < \alpha <1$.
The initial LP problem ($k = 1$) is well-defined, and 
we can give a feasible solution explicitly by
$x_j = 0$ ($j = 1,\ldots, |B|$), $z_{si} = 0$ ($i = 1, \ldots, |V_s|, s = 1, \ldots, |S|$), 
$\overline{\zeta}_s^{\alpha} = 0$ ($\alpha \in \overline{A}_s, s = 1, \ldots, |S|$),
$\underline{\zeta}_s^{\alpha} = 0$ ($\alpha \in \underline{A}_s, s = 1, \ldots, |S|$)
 and $t = \min_{s = 1, \ldots, |S|} \left\{L_s/\underline{P}_s, 
\min_{\alpha \in \underline{A}_s^{\alpha}} \left\{L_s^{\alpha}/\underline{P}_s^{\alpha}, \right\}\right\}$.
Therefore, we know that 
\begin{eqnarray*}
& & \overline{\zeta}_{s}^{\alpha} + \frac{1}{\alpha |V_s| 
 - |\overline{R}_s^{k,\alpha}|} \sum_{i=1,  i \notin \overline{R}_s^{k,\alpha}}^{|V_s|} (z_{si} - \overline{\zeta}_{s}^{\alpha}  )^{+} 
= \overline{\zeta}_{s}^{\alpha} + \frac{1}{\alpha |V_s|} 
\sum_{i=1}^{|V_s|} (z_{si} - \overline{\zeta}_{s}^{\alpha}  )^{+} 
\\ 
& \ge &
 \overline{\zeta}_{s}^{\alpha} + \frac{1}{ |V_s| } 
 \sum_{i=1}^{|V_s|} (z_{si} - \overline{\zeta}_{s}^{\alpha}  )^{+}  
  = 
 \frac{1}{ |V_s|} \sum_{i=1}^{|V_s|} 
\left\{\overline{\zeta}_{s}^{\alpha}  + (z_{si} - \overline{\zeta}_{s}^{\alpha}  )^{+}  \right\} \ge 0.
\end{eqnarray*}
The last inequality holds from an inequality $p + (q-p)^+ \ge 0$ for any $p \in \Real$ and $q \ge 0$.
From (\ref{eq:6c}) and (\ref{New_upper_C}), 
the objective function $t$ has a lower bound 
$ t \ge \max\left\{ \max_{s=1, \ldots, S} \left\{-U_s  / \overline{P}_s\right\},
\max_{\alpha \in \overline{A}_s^{\alpha}} \left\{-U_s^{\alpha}/\overline{P}_s^{\alpha}, \right\} 
\right\}$.
From the duality theorem of linear programming  \cite[etc]{chvatal1983linear},
the initial LP problem has an optimal value $t^{(1)}$.

Next, we assume that the $k$th LP has its optimal value $t^{(k)}$ and optimal solution $x_{i}^{(k)}$ and $z_{si}^{(k)}$,
and we examine the $(k+1)$th LP.
If the number of voxels in $V_s$ such that 
$z_{si}^{(k)} > U_s^{\alpha} + \overline{P}_s^{\alpha} t^{(k)}$ and 
$ i \notin \overline{R}_s^{k, \alpha}$ were no less than 
$\alpha |V_s| - \overline{R}_s^{k, \alpha}$,
we would have
\begin{eqnarray*}
\overline{\zeta}_{s}^{\alpha} + \frac{1}{\alpha |V_s| 
 - |\overline{R}_s^{k,\alpha}|} \sum_{i=1,  i \notin \overline{R}_s^{k,\alpha}}^{|V_s|} (z_{si}^{(k)} - \overline{\zeta}_{s}^{\alpha}  )^{+} 
& > &
\overline{\zeta}_{s}^{\alpha} + \frac{1}{\alpha |V_s| 
 - |\overline{R}_s^{k,\alpha}|} (\alpha |V_s| - |\overline{R}_s^{k, \alpha}|) 
(U_s^{\alpha} + \overline{P}_s^{\alpha} t^{(k)} - \overline{\zeta}_s^{\alpha}) \\
  & = & U_s^{\alpha} + \overline{P}_s^{\alpha}t^{(k)},
\end{eqnarray*}
but this contradicts (\ref{New_upper_C}).
Hence, the number of voxels that will be newly added to $\overline{R}_s^{k+1, \alpha}$ is less than
$\alpha |V_s| - \overline{R}_s^{k, \alpha}$, and this leads to 
\begin{eqnarray}
| \overline{R}_s^{k+1, \alpha} |
 < (\alpha |V_s| - \overline{R}_s^{k, \alpha}) + \overline{R}_s^{k, \alpha} 
= \alpha |V_s|. \label{eq:upperR}
\end{eqnarray}
For the lower DVCs, we also obtain 
$| \underline{R}_s^{k+1, \alpha} | < (1-\alpha) |V_s|$.
Therefore, the denominators in (\ref{New_upper_C}) and (\ref{New_lower_C}) are not zero again, 
and we can use the same discussion as the initial LP problem
to derive that the $(k+1)$th LP has an optimal value $t^{(k+1)}$.
By induction, for any $k \ge 1$, the $k$th LP has its optimal value $t^{(k)}$.

For Part (b),
from the definition of $\overline{R}_s^{k+1, \alpha} =  
\left\{ i \in V_s : z_{si}^{(k)} > U_s^{\alpha} 
+ \overline{P}_s^{\alpha} t^{(k)} \right\}$,
the non-positivity of $t^{(k)}$ indicates that 
$
\{ i \in V_s : z_{si}^{(k)} > U_s^{\alpha} \} \subset 
\overline{R}_s^{k+1, \alpha}$. 
Using the upper bound of the size of $\overline{R}_s^{k+1, \alpha}$ obtained 
in (\ref{eq:upperR}),
we have $| \{ i \in V_s : z_{si}^{(k)} > U_s^{\alpha} \} | \le | \overline{R}_s^{k+1, \alpha}  | \le \alpha |V_s|$ and this means
the solution of $k$th LP problem satisfies the corresponding upper DVCs.
We can also show that the solution with non-positive $t^{(k)}$ satisfies the lower DVCs in a similar way.

Finally, in order to verify the inequality $t^{(k+1)} \le t^{(k)}$ of Part (c), 
we give a feasible solution of the $(k+1)$th LP problem such that $t = t^{(k)}$.
We set $t = t^{(k)}$, $x_j = x_j^{(k)}$, 
$z_{si} = z_{si}^{(k)}$, 
$\underline{\zeta}_s^{\alpha} = L_s^{\alpha} - \underline{P}_s^{\alpha} t^{(k)}$
and  $\overline{\zeta}_s^{\alpha} = U_s^{\alpha} + \overline{P}_s^{\alpha} t^{(k)}$. 
Since these values are derived from the $k$th LP problem,
it is easy to check that these values satisfy the constraints of (\ref{eq:NewProblem}) except
(\ref{New_upper_C}) and (\ref{New_lower_C}).
In (\ref{New_upper_C}), 
the summation $\sum_{i = 1, i \notin \overline{R}_s^{k+1, \alpha}}^{|V_s|} (z_{si} - \overline{\zeta}_s^{\alpha})^+$ is zero due to 
$\overline{R}_s^{k+1, \alpha} =  \left\{ i \in V_s : z_{si}^{(k)} > U_s^{\alpha} + \overline{P}_s^{\alpha} t^{(k)} \right\}$ and 
$\overline{\zeta}_s^{\alpha} = U_s^{\alpha} + \overline{P}_s^{\alpha} t^{(k)}$. Therefore, the left-hand side 
of~(\ref{New_upper_C}) reduces to $U_s^{\alpha} + \overline{P}_s^{\alpha} t^{(k)}$ and this is same as the right-hand side.
Again, we apply a similar step to (\ref{New_lower_C}). 
\hfill \qed


\section{Numerical Experiment}


We used a dataset of the American Association of Physicists in 
Medicine (AAPM) Task Group (TG) 119 report~\cite{ezzell2009imrt}.
 The dataset includes four mock test cases;
a C-shape case, a mock prostate case, a mock head/neck case and a multi target case.
Table~\ref{table:data} is a summary of the dataset.  
For these four cases, the table shows the number of beamlet, the organ names, the number of voxels in the organs and the DVCs.

We compare the proposed method with 
Method-M (the successive LP method of Merrit et al.~\cite{merritt2002successive})
to demonstrate the performance of the proposed method.
We did not include Method-R (the C-VaR method of Romeijn et al.~\cite{romeijn2003novel}) for the comparison,
since we found from preliminary experiments that
the LP problems in the Method-R for all of the four test cases in AAPM TG119 were infeasible.

The dataset of AAPM TG119 is provided as 3D image format called DICOM.
Using the CERR software 4.0~Beta~2~\cite{deasy2003cerr} and MATLAB 2013b,
we transformed the DICOM files into 
the LP problems~(\ref{eq:NewProblem}). 
We ran CERR with its default settings. 
Then, we called CPLEX~12.6.0 to solve the generated LP problems.
Finally, we again utilized CERR to visualize the solutions and checked whether the obtained solutions satisfied the DVCs.
We also used CERR 
to prepare a manageable dataset from the TG119 dataset.
The number of voxels in the PTV organ of the mock head/neck case was
more than 50,000 and this was too large to solve (\ref{eq:NewProblem}) on 16 GB memory space of our computing environment. 
For only this case, therefore, we chose 10,000 voxels randomly from the 50,000 voxels. 
We examined a number of this random selection and we observed this operation did not affect the numerical results so much. 
The computing environment was Windows 8 run on an Intel Core i7-4790 (3.6 GHz, 4 cores) and 16 GB of memory space.

\begin{table}[htbp]
\caption{A summary of the dataset for numerical comparison}\label{table:data}
  \begin{center}
        \begin{tabular}{l|r|lll}
\hline
\hline
	\multicolumn{5}{c}{C-shape (The number of beamlets is 414)} \\
\hline
          organ/tumor & The number of voxels & \multicolumn{3}{c}{DVCs}  \\ \hline
          Outer Target  & 17522  & $L_{\textrm{Outer Target}}^{0.95}$ & $ =$ & $ 50 $ \\ 
                        &        & $U_{\text{Outer Target} }^{0.1}$ & $  = $ & $55 $ \\ 
          Core          &  3087  & $U_{\text{Core} }^{0.1}$ & $ =$ & 
$ 25 $ \\ 
\hline 
\hline
	\multicolumn{5}{c}{Mock Head/Neck (The number of beamlets is 619)} \\
\hline
          organ/tumor   & The number of voxels & \multicolumn{3}{c}{DVCs} \\ \hline
          PTV    & 10000 & $L_{\text{PTV}}^{0.99}$ & $ =$ & $ 46.5$ \\
          &       & $L_{\text{PTV}}^{0.9} $ & $ = $ & $ 50.0  $ \\
          &       & $U_{\text{PTV}}^{0.2}$&$ =$&$ 55  $ \\
          Cord   &  1333 & $U_{\text{Cord}}$&$ =$&$ 40  $ \\
          Lt Parotid & 525 & $U_{\text{Lt Parotid}}^{0.5}$&$ =$&$ 20 $ \\
          Rt Parotid & 740  & $U_{\text{Rt Parotid}}^{0.5}$&$ =$&$ 20 $ \\
\hline 
\hline
	\multicolumn{5}{c}{Mock Prostate (The number of beamlets is 241)} \\
\hline
            organ/tumor & The number voxels & \multicolumn{3}{c}{DVCs} \\ \hline
            ProstatePTV  & 8591  & $L_{\text{ProstatePTV}}^{0.95}$&$ =$&$ 75.6 $ \\ 
                         &       & $U_{\text{ProstatePTV}}^{0.05}$&$  = $&$83  $ \\ 
            Bladder      & 5207  & $U_{\text{Bladder}}^{0.30}$&$ =$&$ 70  $\\ 
                         &       & $U_{\text{Bladder}}^{0.10}$&$ =$&$ 75  $\\ 
            Rectum       & 1830  & $U_{\text{Rectum}}^{0.30}$&$  =$&$ 70  $\\ 
                         &       & $U_{\text{Rectum}}^{0.10} $&$ = $&$ 75  $\\
\hline
\hline
	\multicolumn{5}{c}{Multi Target (The number of beamlets is 601)} \\
            organ/tumor & The number of voxels & \multicolumn{3}{c}{DVCs}  \\ \hline
            Center   & 5143 & $L_{\text{Center}}^{0.99}$&$ =$&$ 50$ \\
                     &      & $U_{\text{Center}}^{0.1} $&$=$&$ 53$ \\
            Superior & 5549 & $L_{\text{Superior}}^{0.99} $&$=$&$ 25$ \\
                     &      & $U_{\text{Superior}}^{0.1} $&$=$&$ 35$ \\
            Inferior & 5529 & $L_{\text{Inferior}}^{0.99} $&$=$&$ 12.5$ \\
                     &      & $U_{\text{Inferior}}^{0.1} $&$= $&$25$ \\
\hline
\hline
          \end{tabular}
        \end{center}
\end{table}

A desirable stopping criterion of the proposed method 
is $t^{(k)} \le 0$, since Theorem~\ref{th:Update} showed 
that the output solution $\x^{(k)}$ satisfies all the DVCs when $t^{(k)} \le 0$.
Due to the intrinsic difficulty arising from combinatorial aspects of the DVCs, 
the number of iterations to attain $t^{(k)} \le 0$ would be prohibitive.
As a practical stopping criterion,
we stop the proposed method when the iteration count $k$ reaches 5
and output $\x^{(5)}$.

For the numerical computation of Method-M (Algorithm~\ref{al:Merritt}), we should describe configurations more specifically.
The description for Method-M in Section~2.2 discussed only one tumor structure and only one healthy structure.
To compute multiple  tumors for Multi Target case, we associated $\tau_s$ for each  tumor structure, 
and we maximized $\sum \tau_s$ over all of the tumor structure.
Next, we set the parameter $\lambda = 10^{-6}$. 
Finally,
we stopped Method-M
when the set $R$ became infeasible on a DVC; more precisely, when $|R_s^{\alpha}|$ exceeded $\alpha |V_s|$.
When the infeasibility was detected at the $k$th iteration, the solution of  Method-M was extracted from the solution 
of $(k-1)$th iteration.

For the execution of the proposed method, we need to specify values
for positive constants $\underline{P}_s, \overline{P}_s,
\underline{P}_s^{\alpha}, \overline{P}_s^{\alpha}$.
We assigned $1$ to these constants for all the cases except the C-shape case.
For the C-shape case, we will explain more details later.

\subsection{Results}

\begin{figure}[htbp]
    \centering
    \includegraphics[width=15cm]{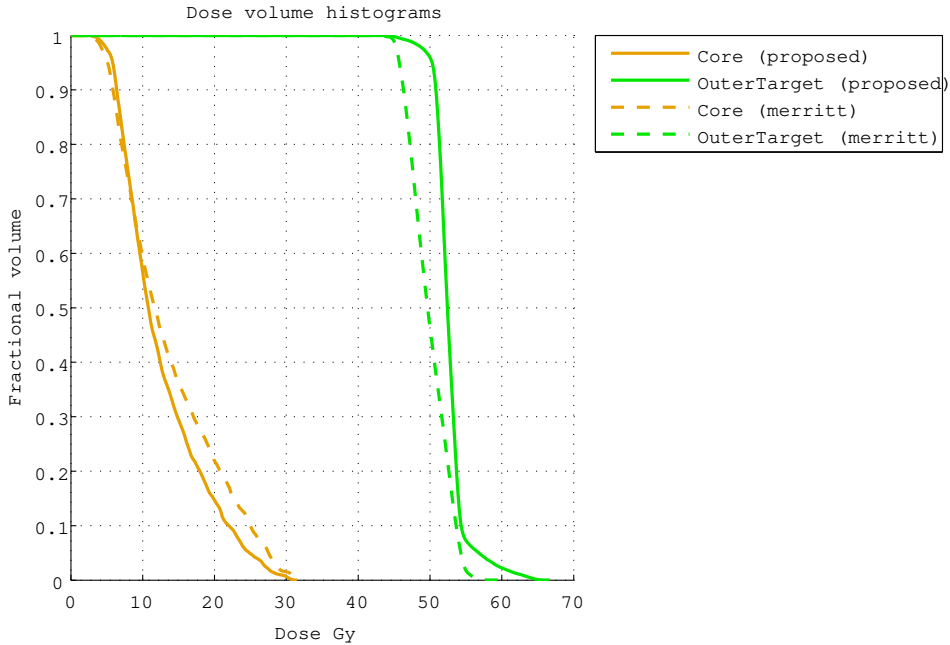}
    \caption{DVH of C-Shape case}
    \label{DVH:C-Shape}
\end{figure}

\begin{figure}[htbp]
    \centering
    \includegraphics[width=15cm]{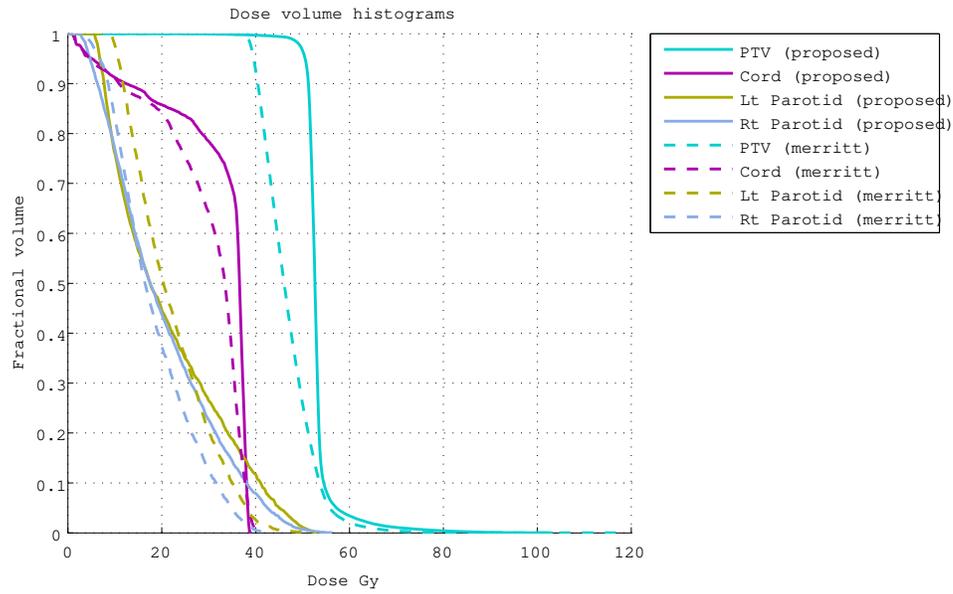}
    \caption{DVH of Head/Neck case}
    \label{DVH:Head}
\end{figure}

\begin{figure}[htbp]
    \centering
    \includegraphics[width=15cm]{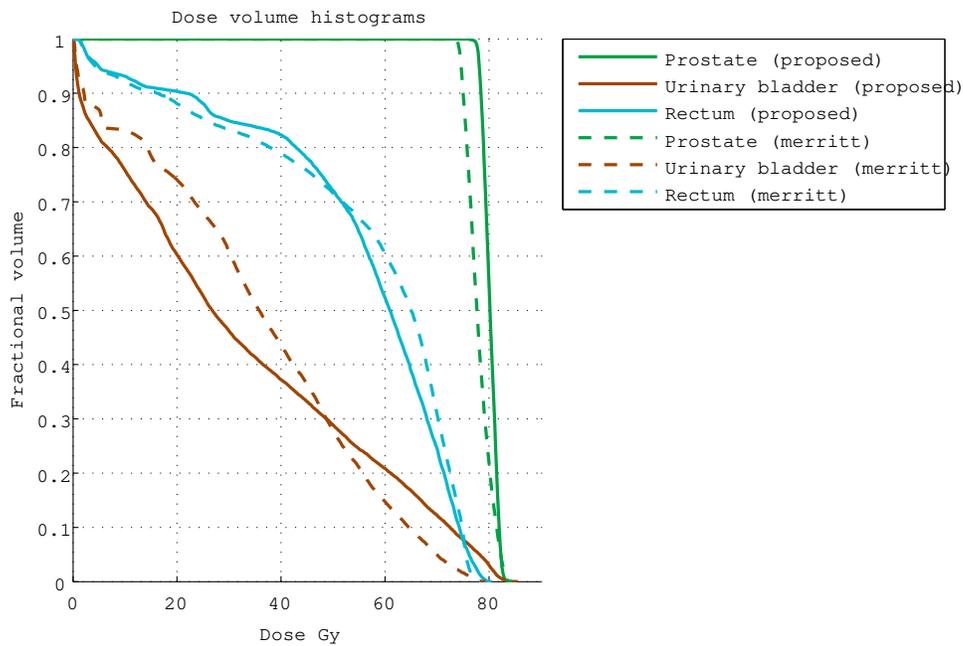}
    \caption{DVH of Prostate case}
    \label{DVH:Prostate}
\end{figure}

\begin{figure}[htbp]
    \centering
    \includegraphics[width=15cm]{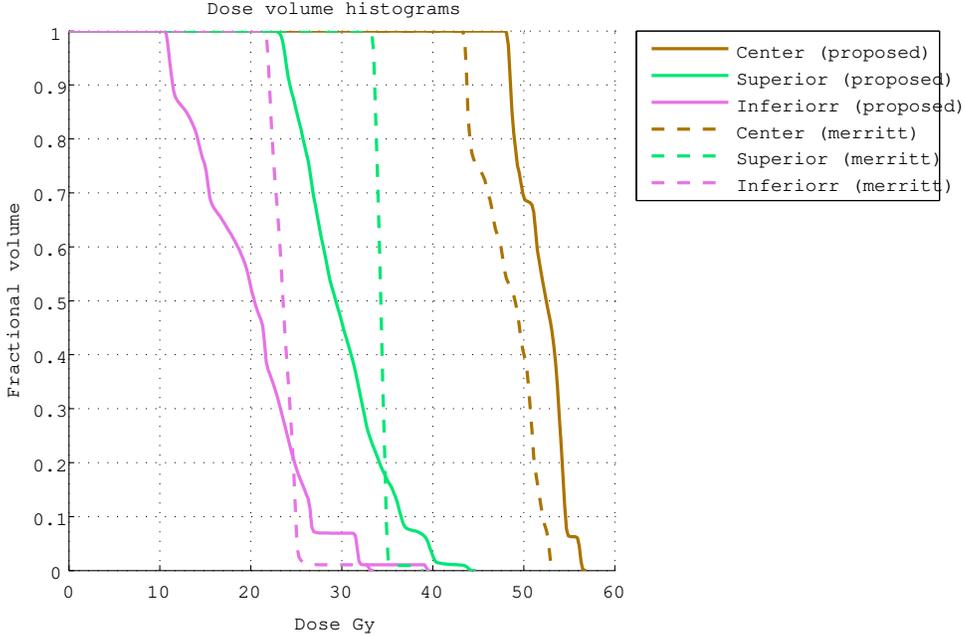}
    \caption{DVH of Multi target case}
    \label{DVH:Multi}
\end{figure}

Table~\ref{table:result} reports whether the obtained solutions satisfies DVCs or not.
In addition, we also provide the DVH figures in Figures~\ref{DVH:C-Shape}, \ref{DVH:Head}, 
\ref{DVH:Prostate}, and \ref{DVH:Multi}.
In these figures, 
the solid lines and the broken lines indicate the results of the proposed method and Method-M, respectively,
and different colors are used to clarify the organs. 

In Table~\ref{table:result}, 
the column ``Proposed''  shows
the evaluation results of the solution obtained by the proposed method
in the viewpoint of DVCs.
For example, the value 50.3 in the row $L_{\textrm{Outer Target}}^{0.95}$ indicates that 
$95\%$ of Outer Target receives at least $50.3$ Gy.
Therefore, this solution satisfies $L_{\textrm{Outer Target}}^{0.95} = 50.0$
and this is indicated by ``Pass.''
The failure of the solutions are indicated by ``Failed'' in the table.
In the same way, the column ``Method-M'' shows the result of Method-M.
The table also reports the numbers of LP problems solved
in each test case. 
The number of the parenthesis in the ``Proposed'' column means that the number of LP problems 
for the proposed method to acquire a feasible solution. 
In C-shape, for example, the proposed method obtained a feasible solution by the third LP problem $(t^{(3)} \le 0)$,
and improved the solution with the successive two LP problems.
(In the Multi Target case, we used $(-)$ to indicate that we failed to obtain a non-negative optimal value.)

 
From the result of Table~\ref{table:result}, 
we observe for the C-Shape case that 
the solution of the proposed method satisfies all DVCs.
Method-M failed in the DVC $L_{\textrm{Outer Target}}^{0.95} = 50.0$.
We also see that the green dashed line 
passes the point $(50,0.5)$ in Figure~\ref{DVH:C-Shape}.
This indicates that in the solution of the Method-M,
a half of the voxels in OuterTarget 
receives less than 50Gy. 
We can further acquire similar observations on the Head/Neck case and the Prostate case.

For the Multi Target case, however,
Table~\ref{table:result} shows that
the solution of the proposed method
fails to satisfy all DVCs.
This means that the solution is far from 
a favorable solution.
In addition, the shapes of the proposed method in Figure ~\ref{DVH:Multi} are gradual slopes.
In contrast,
Method-M has a narrow width and outputs a better solution than the proposed method for only the Multi Target case.

Since the proposed method outputs solutions
that matches the DVCs more adequately than Method-M for three problems out of the four cases, 
the proposed method has a tendency to output favorable solutions than Method-M.
Furthermore, Table~\ref{table:result} shows that the numbers of iterations in the proposed method were
at most 5. 
This implies that the proposed  is better than Method-M from the viewpoint of computation cost.

\begin{table}[htbp]
 \caption{Numerical results for the four test cases}
    \label{table:result}
  \begin{center}
    \begin{tabular}{l|r|l|l}
\hline
\hline
\multicolumn{4}{c}{C-Shape} \\
\hline
      organ/tumor  & \multicolumn{1}{|c|}{DVCs} & \multicolumn{1}{|c|}{Proposed} & \multicolumn{1}{|c}{Method-M} \\
\hline
      Outer Target & $L_{\textrm{Outer Target}}^{0.95} = 50.0 $   & 50.3 (Pass) & 45.7 (Fail) \\
                   & $U_{\textrm{Outer Target}}^{0.10} = 55.0 $   & 54.3  (Pass) & 53.6 (Pass)   \\ 
      Core         & $U_{\textrm{Core}}^{0.10} = 25.0$  & 22.0 (Pass) & 25.0 (Pass) \\
\hline
\multicolumn{2}{l|}{the number of iterations} & 5(3) & 14 \\
\hline
\hline
\multicolumn{4}{c}{Mock Head/Neck} \\
\hline
      organ/tumor  & \multicolumn{1}{|c|}{DVCs} & \multicolumn{1}{|c|}{Proposed} & \multicolumn{1}{|c}{Method-M} \\
\hline
      PTV & $L_{\textrm{PTV}}^{0.99} = 46.5 $   & 47.6 (Pass) & 38.4 (Fail) \\
                   & $L_{\textrm{PTV}}^{0.90} = 50.0 $   & 51.2  (Pass) & 40.5 (Fail)   \\ 
                   & $U_{\textrm{PTV}}^{0.20} = 55.0 $   & 53.5  (Pass) & 51.2 (Pass)   \\ 
      Core         & $U_{\textrm{Core}} = 40.0$  & 39.1 (Pass) & 41.3 (Fail) \\
      Lt Parotid         & $U_{\textrm{Lt Parotid}}^{0.50} = 20.0$  & 17.5 (Pass) & 20.4 (Fail) \\
      Rt Parotid         & $U_{\textrm{Rt Parotid}}^{0.50} = 20.0$  & 17.5 (Pass) & 16.4 (Pass) \\
\hline
\multicolumn{2}{l|}{the number of iterations} & 5(3) & 27 \\
\hline
\hline
\multicolumn{4}{c}{Prostate} \\
\hline
      organ/tumor  & \multicolumn{1}{|c|}{DVCs} & \multicolumn{1}{|c|}{Proposed} & \multicolumn{1}{|c}{Method-M} \\
\hline
      Prostate PTV & $L_{\textrm{Prostate PTV}}^{0.95} = 75.6 $   & 78.1 (Pass) & 74.7 (Fail) \\
                   & $U_{\textrm{Prostate PTV}}^{0.05} = 83.0 $   & 82.3  (Pass) & 82.2 (Pass)   \\ 
      Bladder         & $U_{\textrm{Bladder}}^{0.30} = 70.0$  & 48.8 (Pass) & 48.7 (Pass) \\
               & $U_{\textrm{Bladder}}^{0.10} = 75.0$  & 72.5 (Pass) & 64.5 (Pass) \\
      Rectum         & $U_{\textrm{Rectum}}^{0.30} = 70.0$  & 68.1 (Pass) & 70.2 (Pass) \\
           & $U_{\textrm{Rectum}}^{0.10} = 75.0$  & 74.3 (Pass) & 74.8 (Pass) \\
\hline
\multicolumn{2}{l|}{the number of iterations} & 5(2) & 22 \\
\hline
\hline
\multicolumn{4}{c}{Multi Target} \\
\hline
      organ/tumor  & \multicolumn{1}{|c|}{DVCs} & \multicolumn{1}{|c|}{Proposed} & \multicolumn{1}{|c}{Method-M} \\
\hline
      Center & $L_{\textrm{Center}}^{0.99} = 50.0 $   & 48.2 (Fail) & 43.5 (Fail) \\
                   & $U_{\textrm{Center}}^{0.10} = 53.0 $   & 54.6 (Fail) & 52.1 (Pass)   \\ 
      Inferior         & $L_{\textrm{Inferior}}^{0.99} = 12.5$  & 10.8 (Fail) & 21.7 (Pass) \\
               & $U_{\textrm{Inferior}}^{0.10} = 25.0$  & 26.6 (Fail) & 24.9 (Pass) \\
      Superior        & $L_{\textrm{Superior}}^{0.99} = 25.0$  & 23.3 (Fail) & 33.3 (Pass) \\
           & $U_{\textrm{Superior}}^{0.10} = 35.0$  & 36.5 (Fail) & 34.9 (Pass) \\
\hline
\multicolumn{2}{l|}{the number of iterations} & 5(-) & 15 \\
\hline
    \end{tabular}
\newline
  \end{center}
\end{table}

\section{Discussions}

In this paper we propose a new iterative algorithm for the FMO problems.
Here, we discuss several aspects of the proposed method.

\subsubsection*{The result on Multi Target}
The results on the Multi Target case
imply an weakness of the proposed method.
We have two reasons that the proposed method did not work properly for this test case.
First, the concept of our method is to get a better solution by adjusting unsatisfied inequalities in
the LP problem~(\ref{eq:NewProblem}).
Hence, the solution in a next iteration often improves the satisfiability of the unsatisfied inequalities.
At the same time, however, this may bring a negative effect on the inequalities that are already satisfied.
In particular, if there are hard constraints (with $L_s$ and $U_s$) and soft constraints (with $L_s^{\alpha}$ and $U_s^{\alpha}$)
together in one problem, 
the soft constraints are influenced by the hard constraints and
they sometime become infeasible even when there are feasible solutions that fulfill these soft constraints.
Second, Ezzel et al.~\cite{ezzell2009imrt} pointed out that the Multi Target case is a difficult test case.
They confirmed from statistics that for cases like the Multi Target case,
even treatment planners often generate a solution that do not satisfy DVCs.
For such cases, we should employ a different criterion instead of pursuing DVCs.

\subsubsection*{Parameters for Method-M}
There would be a possibility that a careful selection on the parameters for Method-M might improve the quality of the solution or the running time.
In the numerical experiments, we examined Method-M changing the parameters 
and we reported the best results of Method-M from the different parameters,
so further improvements only by the parameter selection are not so promising.

For Method-M, a reduction of the computation time is also a daunting task.
The test cases we used have a few DVC. 
To solve one LP problem, the proposed method consumed 
 2-3 times computation time of Method-M.
However, the proposed method acquired a feasible solution 
in two or three iterations.
To compete the proposed method, therefore,
Method-M should complete its algorithm
with at most ten iterations, but Method-M required at least 14 iterations as shown in  Table~\ref{table:result}.
We remark that the number of intermediate variables in the LP problems of the proposed method 
depends on the number of DVCs, while that of Method-M is independent from the number of DVCs.
Therefore, Method-M may perform well in a test case with a large number of DVCs.

\subsubsection*{An extension of our approach for precise volumes}
Another issue from a different viewpoint is the volume of each voxel.
In this paper, we assumed that an organ was divided into 
voxels of the same rectangular shape, thus all the voxels had the same volume.
However, voxels at the boundary of an actual structure 
may partially contain exterior of the structure.
Therefore, there may be a difference between the total volume of voxels in the structure
and the actual volume,
and  our approach would output a solution with serious deviations.

Our approach can be extended to 
handle the precise volume of each voxel by the steps below.
For the $i$th voxel of the structure $s$, 
$c_{si}$ is used to denote the volume of a part of  $s$ that is covered by 
the $i$th voxel, in other words, $c_{si}$ is the volume of the intersection 
of $s$ and the $i$th voxel.
Then, a constant value $C_s := \sum_{i=1}^{v_s} c_{si}$ 
denotes the actual volume of  $s$.
We also define a new set for outliers  $\overline{C}_{s}^{k, \alpha} := \sum_{i \in \overline{R}_{s}^{k, \alpha}} c_{si}$ for an upper DVC. In addition, $\underline{C}_{s}^{k, \alpha}$ is introduced for an lower DVC．

We replace the $k$th LP problem in Algorithm~\ref{al:proposed} with the following LP problem:
\begin{subequations}
\label{Extended}
\begin{align} 
  &&&&&\textrm{minimize}    && t & &&&&& \\
  &&&&& \textrm{subject to} && \sum_{j=1}^{N_b} D_{sij} x_j = z_{si} & i \in V_s;  s \in S \\
  &&&&&                     && L_s -\underline{P}_s t \leq z_{si}\leq U_s + \overline{P}_s t &  i \in V_s; s \in S \label{New_all}\\
  &&&&&                     && \underline{\zeta}_{s}^{\alpha} - \frac{1}{(1-\alpha) C_s - \underline{C}_{ks}^{\alpha} } \sum_{i \notin \underline{R}_{s}^{k, \alpha}} c_{si} ( \underline{\zeta}_{s}^{\alpha} - z_{si} ) ^{+}  \geq L_s^{\alpha} - \underline{P}_s^{\alpha}t  & \alpha \in \underline{A}_s ; s \in S \label{Extend_lower}\\
  &&&&&                     && \overline{\zeta}_{s}^{\alpha} + \frac{1}{\alpha C_s -  \overline{C}_{s}^{k, \alpha}  } \sum_{i \notin \overline{R}_{s}^{k, \alpha}} c_{si}(z_{si} - \overline{\zeta}_{s}^{\alpha}  )^{+}  \leq U_s^{\alpha} + \overline{P}_s^{\alpha}t & \alpha \in \overline{A}_s ; s \in S \label{Extend_upper}\\
  &&&&&                     && x_j  \geq 0   & j \in B &&&&& \\
  &&&&&                     && z_{si}  \geq 0  & i \in V_s; s \in S    &&&&& \\
  &&&&&                    &&  \underline{\zeta}_{s}^{\alpha},  \overline{\zeta}_{s}^{\alpha},t : \textrm{free variables} &&
\end{align}
\end{subequations}

A main difference between the original proposed method and this extended method lies in (\ref{New_lower_C}) and (\ref{Extend_lower}). In the original proposed method, we use the number of voxels to represent fractional volume of a structure. 
However, in the extended method, we use their actual volumes of voxels. Therefore, we include $c_{si}$ in the summation of (\ref{Extend_lower}). 
This LP problem satisfies the same property as Theorem~\ref{th:Update}.
In particular, we can find a solution that satisfies all the DVCs
when the optimal value $t^{(k)}$ of the $k$th LP problem is non-positive.
Therefore, we can naturally extend the proposed method to handle the actual volumes.


\section{Conclusion and Future Directions}

In this paper, we proposed a new method for FMO problems
using the C-VaR type constraints of Romeijn et al.
and the outliers of Merrit et al. 
The proposed method has favorable mathematical properties as discussed in Theorem~\ref{th:Update}. In particular, when the optimal value of the LP problems is non-positive, 
its optimal solution satisfies all the DVCs.
From the numerical experiments, we verified that 
the proposed method was effective for the mock case.
Particularly, the proposed method found feasible solutions
for the test cases whose feasibilities were not detected by the LP problems of Romeijn et al.
Furthermore, our approach obtained these feasible solutions
within a shorter computation time than the approach of Merrit et al.

Further studies should include the removal of the two weaknesses in the proposed method.
The first one is to get a better solution for the Multi-Target  case and similar cases.
Though these cases are very hard as discussed in Section~5,
deep investigation the effect of the parameters 
involved in the LP problems could improve the situation. In particular, we need a discussion on adaptive selections of 
$\underline{P}_s^{\alpha}$ and $\overline{P}_s^{\alpha}$.
The other issue is a strong dependence of the size of LP problems 
on the number of DVCs. 
The size of the LP problems in the proposed method
grows rapidly, when the number of DVCs increases.
As a result, the size of solvable FMO problems are limited.
We should reduce the number of variables involved in the LP problems
by detecting redundant variables or inactive constraints.

%


\end{document}